\documentclass{amsart}
\usepackage{graphicx}
\graphicspath{ {.} }
\setkeys{Gin}{width=0.8\linewidth,keepaspectratio=true}
\usepackage[latin1]{inputenc}
\usepackage{natbib}
\bibpunct{(}{)}{;}{a}{,}{;}
\usepackage{xcolor}
\usepackage{url}
\usepackage{epsfig}
\usepackage{algorithm}
\usepackage{algorithmic}
\usepackage[algo2e]{algorithm2e}
\usepackage{longtable}

\usepackage{rotating}
\usepackage{subfigure}
\usepackage{lscape}
\usepackage{multirow}
\usepackage{adjustbox}
\usepackage{comment}
\usepackage{setspace}

\begin{document}

\title[Amazon Locker capacity management]{Amazon Locker capacity management}

\author[S. Sethuraman]{S. Sethuraman}
\address[Samyukta Sethuraman, \textcolor{black}{ORCID 0000-0003-2803-6550}]{
Sponsored Products Advertising, Amazon,  Palo Alto, CA}
\email[S. Sethuraman]{samyukts@amazon.com}

\author[A. Bansal]{A. Bansal}
\address[Ankur Bansal, \textcolor{black}{ORCID 0000-0002-1842-5125}]{Kotak 
Mahindra Bank, Gurugram, Haryana, India}
\email[A. Bansal]{ank.bnsl@gmail.com}

\author[S. Mardan]{S. Mardan}
\address[Setareh Mardan, \textcolor{black}{ORCID 0000-0003-1515-2093}]{Retail Pricing, Pricing Research, 
Amazon, Seattle, WA}
\email[S. Mardan]{msetareh@amazon.com}

\author[M.G.C. Resende]{M.G.C. Resende}
\address[Mauricio G.C. Resende, \textcolor{black}{ORCID 0000-0001-7462-6207}]{Industrial \& Systems Engineering, University of Washington, 
Seattle, WA}
\email[M.G.C. Resende]{mgcr@uw.edu}

\author[T.L. Jacobs]{T.L. Jacobs}
\address[Timothy L. Jacobs, \textcolor{black}{ORCID 0000-0001-9403-7100}]{ATS Science \& Engineering, Amazon Transportation Services,
Amazon, Bellevue, WA}
\email[T.L. Jacobs]{timojaco@amazon.com}

\begin{abstract}
Amazon Locker is a self-service delivery or pickup location where
customers can pick up packages and drop off returns. A basic 
first-come-first-served policy for accepting package delivery requests
to lockers results in lockers becoming full with standard shipping
speed (3--5 day shipping) packages, and leaving no space left for
\textcolor{black}{expedited} packages which are mostly Next-Day 
or Two-Day shipping. 
This paper proposes a solution to the problem of determining
how much locker capacity to reserve for different ship-option packages.
Yield management is a much researched field with popular
applications in the airline, car rental, and hotel industries.
However, Amazon Locker poses a unique challenge in this field since
the 
\textcolor{black}{number of days a package will wait in a locker (package
dwell time)} is, in general, unknown. 
The proposed solution combines
machine learning techniques to predict locker demand and package
dwell time, and linear programming to maximize throughput in lockers.
The decision variables from this optimization provide optimal
capacity reservation values for different ship options. This resulted
in a year-over-year increase of 9\% in Locker throughput worldwide
during holiday season of 2018, impacting millions of customers.
\end{abstract}
\date{December 11, 2023}
\thanks{To appear in INFORMS Journal on Applied Analytics.}
\maketitle

\section{Introduction}

Amazon Locker is a self-service \textcolor{black}{parcel locker} that accepts
package deliveries and customer returns. 
\textcolor{black}{They} are conveniently
located in highly-visited locations such as shopping centers,
business parks, gas stations and convenience stores. 
On the Amazon.com order
placement page, customers have the option of getting a package
delivered to a locker of their choice instead of their home, at no
additional cost. 
\textcolor{black}{
The selection of lockers presented to a customer includes those
with available space and situated close to the customer's residence.
After a locker is allocated for the order, the package is delivered
there, and the customer receives a notification along with a code
to access a specific bin within the locker. The notification also
indicates the timeframe within which the package must be picked up.
If the package is not retrieved within this period, a return is
initiated and a credit is issued to the customer.
The number of days the package remains in the locker is
the random variable \textit{dwell time}.}
\textcolor{black}{Lockers}  provide unique advantages such as protection
from package theft, and the option of pickup or drop off at a time
and location convenient to the customer. 
Lockers have become
increasingly popular among Amazon customers due to these unique
features.



\begin{figure}[hp]
\centering
\begin{tabular}{c}
\hline
\multicolumn{1}{|c|}{\includegraphics[scale=0.5]{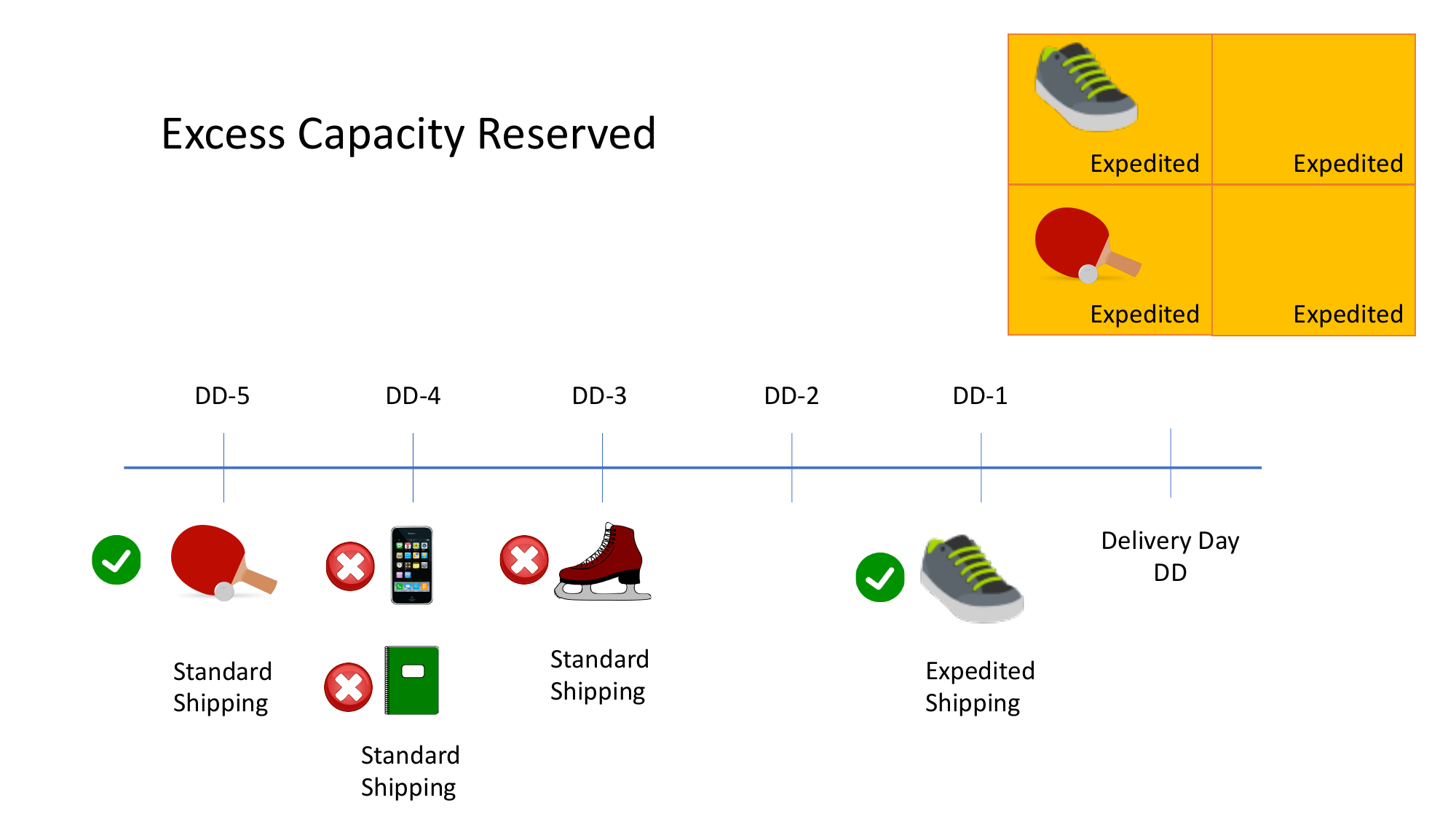}} \\
\hline
\vspace*{5pt} \\
\hline
\multicolumn{1}{|c|}{\includegraphics[scale=0.5]{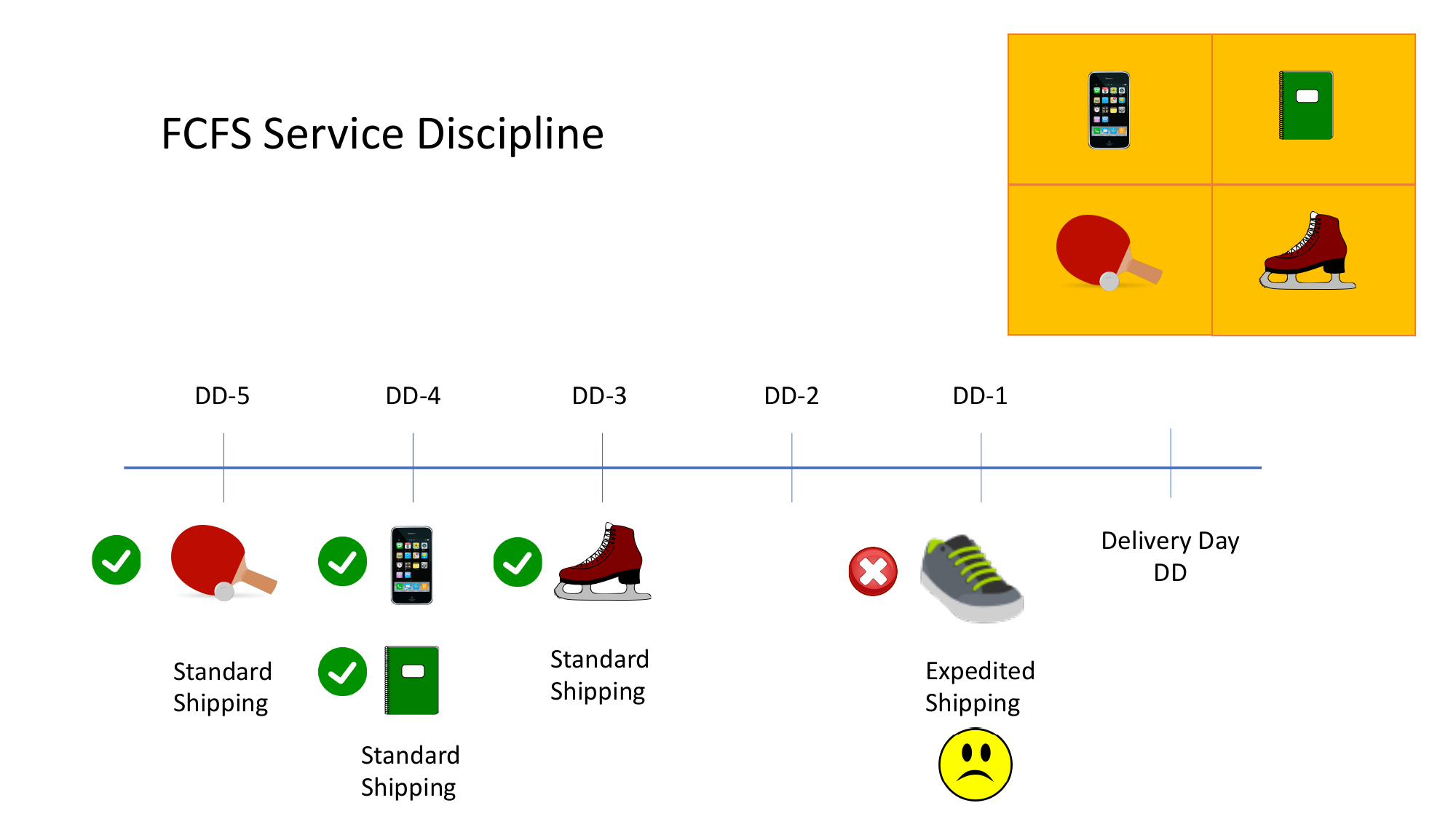}} \\
\hline
\vspace*{5pt} \\
\hline
\multicolumn{1}{|c|}{\includegraphics[scale=0.5]{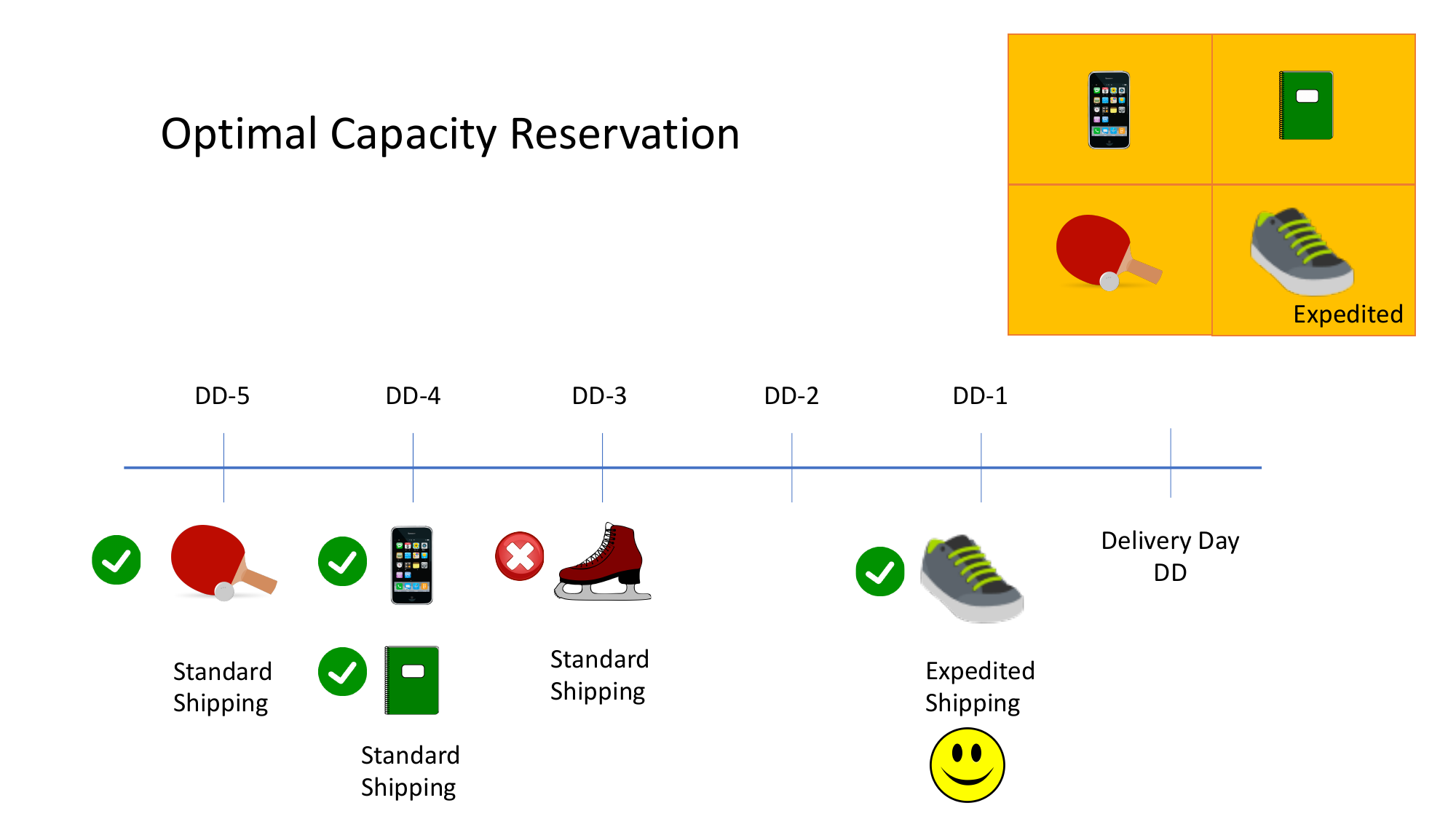}} \\
\hline
\end{tabular}
\caption{Three distinct Locker Capacity Management disciplines.  On top 
excess capacity is reserved for \textcolor{black}{expedited}  shipments \textcolor{black}{and locker is underutilized}.  In the middle, the 
first-come-first-served (FCFS) rule is implemented \textcolor{black}{and customer ordering expedited cannot ship to locker}. 
An optimal capacity reservation scheme is
shown in the bottom.}
\label{f_example-excess}
\end{figure}

With increased demand, capacity management has become an essential
component of the operations of Amazon Locker. 
Locker capacity management performs
two major functions:

\begin{itemize}
\item \textit{Delivery Request Evaluation}: The time lag between when
a customer places an order to when it is delivered can vary between
a few hours and
a few days, depending on the ship option chosen.
Customers have the option
of choosing one out of multiple ship options which determines the
shipping speed of the package. For example, in the US, order ship
options can
be Same-Day, Next-Day, Two-day, or Standard (3-5 days).
Customers can also choose to return their orders to a
locker, which can be considered as a special type of ship option.
When a customer requests an order to be delivered to a locker, the
capacity management system predicts (i) when the package will be
delivered, and (ii) the occupancy of the locker at the time when
the package will be delivered. 
Based on these predictions, the
request is accepted if capacity is estimated to be available, or
rejected otherwise.
\item \textit{Ship Option Capacity Reservations}: 
Typically, if the time-line of order requests for delivery to a
locker on a particular day is concerned, slower ship option orders
are placed before faster ship option orders. 
Therefore, a
first-come-first-served service discipline for accepting delivery
requests to lockers results in the locker running out of space and
the rejection of faster ship option packages. 
As we will see next, this is not desirable if the goal is 
throughput maximization
\textcolor{black}{
since the average number of days a package stays in the locker (the
package dwell time) is usually lower for faster ship option packages.
Therefore, accepting a higher number of faster ship option packages results in
higher locker throughput and subsequently, higher number of satisfied
customers.}

\begin{figure}[ht]
\centering
\includegraphics[scale=1]{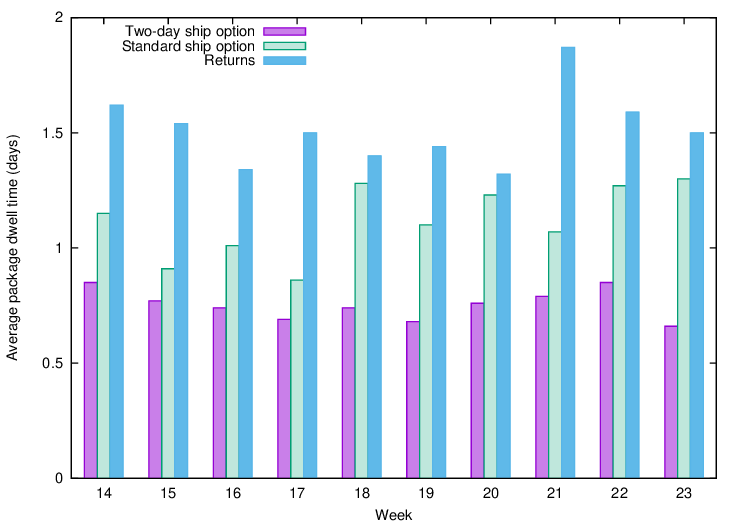}
\caption{
Average package dwell time (in days) for Amazon Locker locations in zip-code 98109 (Seattle, Washington) during ten weeks in 2018.}
\label{f_package_dwell_time}
\end{figure}

Figure~\ref{f_package_dwell_time} shows average
dwell time in days for Two-Day, Standard and Customer Return packages
over a 10-week period in 2018 for lockers in the zip-code 98109 (Seattle, Washington).
\textcolor{black}{
The term \textit{dwell time} is the duration (in days) that a package 
remains in the locker until it is picked up.
}
In this figure, dwell time in days
for packages picked up on the same day as delivery is considered to be zero. 
The data show that dwell time of \textcolor{black}{two-day ship option packages}
is consistently lower
than that of standard, which in turn is lower than the dwell time
of customer returns.
Therefore, it is necessary to reserve space for faster ship option
packages. 
However, reserving more than the required space for faster
ship options 
unnecessarily
results in rejection of slower ship option requests
and therefore results in wasted locker space.
Figure~\ref{f_example-excess} shows three examples of Locker capacity
management disciplines.
\textcolor{black}{
In the time axes of this figure, DD
is the delivery day, i.e. the day the package is placed in the locker,
while DD\hspace*{1pt}-\textit{i} 
indicates $i$ days before the delivery date.
In the three examples, the yellow rectangle represents a locker with 
four slots.
The timeline shows when orders were made.}
In the top example, excessive capacity is reserved for possible 
\textcolor{black}{next-day or two-day ship option
customers.}  
In this example, only the standard shipped ping-pong racket order,
\textcolor{black}{made
five days prior to the delivery date,} and
the \textcolor{black}{next-day} shoe \textcolor{black}{order}
are allocated a slot in the Locker while the mobile phone,
the book, and the ice skate are rejected.
\textcolor{black}{
In this example, two locker slots are left unused.}
In the other extreme, the example in the middle of the figure, a 
first-come-first-served discipline is shown
where the orders are assigned a slot
\textcolor{black}{ as the orders
are placed}.  
In this discipline \textcolor{black}{expedited} orders are usually rejected because the
locker is usually full with standard shipping orders when a new 
\textcolor{black}{expedited} order 
is made.  
These two extreme cases motivate our solution, shown in the bottom of the 
figure.
In an optimal capacity reservation approach 
\textcolor{black}{next-day or two-day} orders are accepted 
without the need to use the excessive capacity 
reservation that forces standard shipping 
orders to be rejected.
This leads to great customer satisfaction since both 
\textcolor{black}{expedited} and standard
shipping orders are accepted by the Locker.
\end{itemize}

This paper describes the algorithm developed by us at Amazon to enable
the second function of locker capacity management described
above.
We answer the following question for each locker: How much capacity should be
reserved for each ship option on each day so
as to maximize the throughput of the locker?
\textcolor{black}{
In the context of lockers, especially delivery or parcel lockers,
\textit{throughput} generally refers to the number of packages (or
transactions) that a locker can handle or process within a given
period of time. Throughput can be used as a metric to gauge the
efficiency and capacity of a locker system.}

At a high level, this problem falls into the widely researched field
of yield management. 
Yield management principles have been applied
widely in travel industries such as airlines and hotels. 
One of the first publications on applied airline
yield management is the seminal paper by
\citet{Lit72a,littlewood} which
introduces Littlewood's rule for assigning capacity
reservations to two fare classes in single leg flights.
This was
extended to multi-leg flights with multiple fare classes 
by \citet{wang1983optimum} in what
later
became known as Expected Marginal Seat Revenue (EMSR).
Optimization and estimation problems 
in airline
yield management are addressed
in a PhD thesis by \citet{McG89a}.
The 1991 Franz
Edelman Award was won by American Airlines. The award-winning work
is published in \citet{smith1992yield}, where a complete implementation
of an end-to-end yield management system at American Airlines is
presented. 
\citet{jacobs2008incorporating} present scalable
yield management models that are currently used by modern-day
airlines. 
\textcolor{black}{
Yield management at Hub Group,
a North American freight rail transportation company, 
is discussed in \citet{Gor10a}. The paper describes
an integrated decision support system to enhance yield 
management and container allocation. 
In its first year of use by Hub Group in 2008, the system increased revenue 
per load by 3\%, container velocity by 5\%, 
and generated \$11 million in cost savings, yielding a 22-fold return 
on the initial investment.}
\citet{McGRyz99a} review over 190 references on revenue management prior to
1999.
\citet{WeaKim03a}
used real data to test a variety of forecasting methods
in a hotel revenue management system.
\citet{BitCal03a}
survey research on dynamic pricing policies and their
relation to revenue management up to 2002.
\textcolor{black}{
\citet{GuiRuiDelMacMacPinPro19a} 
describe Opticar, a decision support system that uses advanced algorithms 
to forecast demand, optimize revenue, and manage Europcar's car rental 
fleet capacity for up to six months in advance. The system, adopted by managers in nine countries, serves as a foundation for daily operations, 
aiding in pricing and vehicle management.
\citet{TalRuzKarVul09a} provide an introduction to revenue management from 
the perspective of simulation.
\citet{GopRan10a}
present
a linear programming model to optimally
allocate train capacity among different travel segments on an Indian
Railways route. Applying this model to 17 trains resulted in
significant increases in revenue, load factors, and passengers
carried.
Recognizing the distinct needs of the cruise industry,
\citet{BecEtAlii21a} describe 
the Yield Optimization and Demand Analytics (YODA) system,
a revenue management systems at
Carnival Corporation. 
This system employs a unique model and machine learning to 
determine cruise pricing and inventory allocation, 
driving a 1.5\% to 2.5\% rise in net ticket revenue.
}
\citet{CloJacGel15a}
address the network revenue management problem under a 
MultiNomial Logit
discrete choice model with applications in airline revenue management.
\citet{FelTop17a}
study revenue management problems where customers select from
a group of offered products according to the Markov chain choice model.
\citet{FerSimWan18a}
consider a price-based network revenue management problem 
where the goal of
a retailer is to maximize revenue
over a finite selling season
from multiple products 
with limited inventory.
\citet{WanLiuCheRenMenHu18a} explore revenue management in
online advertisement.
\citet{KunTal19a}
propose a new Lagrangian relaxation method for discrete-choice 
network revenue management
based on an extended set of multipliers.
\citet{ShiLogThoWei19a} consider
the seat inventory control and
over-booking problem of airline revenue management,
formulated as a Markov Decision
Process and solved with Deep Reinforcement
Learning to find a policy that
maximizes the revenue of each flight.
\textcolor{black}{
\citet{PanDutJos17a} reviews revenue-management research in the 
broadcasting and online advertisement industries, focusing on 
strategies and techniques to maximize advertising revenue. 
It also highlights mobile advertising as an emerging area in revenue 
management and identifies potential gaps for future research.
}
The book by
\cite{phillips2005pricing}
provides a comprehensive review of practical yield management in
industry.

\textcolor{black}{
Parcel lockers have been studied in the literature since the mid-2010s.
A literature review of parcel lockers can found in \citet{DuiWieAreAms20a}.
In that paper,
four literature focus areas are identified: 1) Parcel locker use from the customer's
perspective; 2) Location of parcel lockers; 3) Cost of parcel lockers; 
and 4) Environmental economics impact.
A case study of PostNL, the market segment leader (70\%)
in The Netherlands, illustrates the four focus areas.
\citet{IwaKijLem16a} concentrates on the analysis of evaluating the
usability and efficiency 
of parcel lockers using the Polish InPost Company system as a case study. 
\citet{OreRavSad19a} introduces a logistic model for small parcel delivery 
using multiple service points (SPs) and showcases solution methods 
like the savings heuristic, the petal method, and tabu search. 
The model, demonstrating cost and time efficiencies, especially when 
recipients are flexible with delivery locations, is validated through 
numerical and simulation studies against traditional methods.
}

\textcolor{black}{
\citet{RohGen20a}
explores various locker station delivery concepts, addresses related
challenges and decision problems, and outlines potential research
directions in operations research (OR) to foster advancements in
this burgeoning domain.
In Section~3.3. of their monograph, \textit{Locker assignment and scheduling},
they conclude with the recommendation 
``More research is therefore needed, focusing on 
the capacity management of the locker stations at the operational level, while 
taking into account uncertainties and the dynamic operating environment.''
Our paper is a first step in this direction.
}

Yield management in Amazon Locker offers some unique challenges,
not common in the travel industry. 
The major challenge arises
from unknown package dwell times, i.e., the time a package
stays in the locker (which is known at the time of booking for
airlines, \textcolor{black}{car rentals}, and hotels). 
Furthermore, the sparsity of data
\textcolor{black}{(due to the fact that lockers are deployed in areas with 
different demographics, where pickup and ordering characterists vary in 
volume and timing)}, and the
super-linear growth rate of Amazon Locker locations 
offer additional challenges to
traditional assumptions on demand distribution. 
These factors forge
the need for disruptive technologies for \textcolor{black}{locker} capacity
management.

Apart from unique modeling challenges, the need for such
disruptive
technologies is reinforced by the impact of Amazon Locker capacity management
on customers. 
The capacity management system previously in place
\textcolor{black}{at Amazon (the \textit{legacy} system)} 
resulted in many locker requests being rejected even when capacity
was available in the locker. 
Such rejected requests are called
\textit{unjustified rejections}. 
These rejections happen because of poor locker demand
and dwell time forecasts, resulting in poor ship option capacity
reservations. 
In 2017, twenty percent of Locker requests were rejected.
Out of these $20\%$, between \textcolor{black}{$18\%$} to $35\%$ were unjustified rejections.
Forty percent of these unjustified rejections
were attributed to poor ship option capacity reservations.
Customer surveys showed that these rejections were the biggest pain point
for Amazon Locker customers.

As the Amazon Locker business geared up for phenomenal growth, it was
important to invest into efficient capacity management to
ensure the best possible utilization of resources and guarantee great
customer experience.

The remainder of the paper is organized as follows.
In Section~\ref{s_current_practice} we review the methodology previously
applied at Amazon for locker capacity management.
Section~\ref{s_model} focuses on the newly proposed Locker Capacity
Management model comprised of a locker demand forecasting module, a package
dwell time probability estimation module, and a module for ship option
capacity reservation optimization.
In Section~\ref{s_simulation} we describe a data-driven
simulation system used to evaluate 
our proposed model and
determine the impact on the locker capacity management system of any 
changes made to inputs, 
including changes in ship option capacity reservations.
In Section~\ref{s_results} we report experimental results 
from the implementation of the Locker Capacity Management Model 
both in the simulation system and in real-world production.
Finally, conclusions are drawn in Section~\ref{s_conclusion}.

\section{Legacy Practice} 
\label{s_current_practice}

The legacy methodology for assigning capacity reservations for
faster ship options was based on the following algorithm.
For each locker, in each week, evaluate the proportion of different
ship option packages delivered to homes in the previous year, in
that week, in the zip-code of the locker. 
Locker capacity was then
assigned as equal to these proportions, normalized by locker capacity.
We refer to this algorithm in the remainder of the paper as the
\textit{proportion rule}.

As a hypothetical example of the proportion rule, 
consider the locker Ruby (located in
Amazon's Ruby building in zip-code 98109 in Seattle, Washington). 
Assume that
there are only two ship options to consider: Two-Day and Standard.
If the number of packages delivered to homes in the zip-code 98109,
in week $26$ of $2017$, was $3000$ Two-Day packages and $1000$
Standard packages, and if the locker capacity is $100$ slots, then 
in week $26$ of $2018$ (one year later),
$75$ $(= 100 \times 3000/(3000+1000))$
slots would be reserved for Two-Day packages and 
$25$ $(= 100 \times 1000/(3000+1000))$ 
for Standard packages.

\begin{figure}[t]
\centering
\includegraphics[scale=0.5]{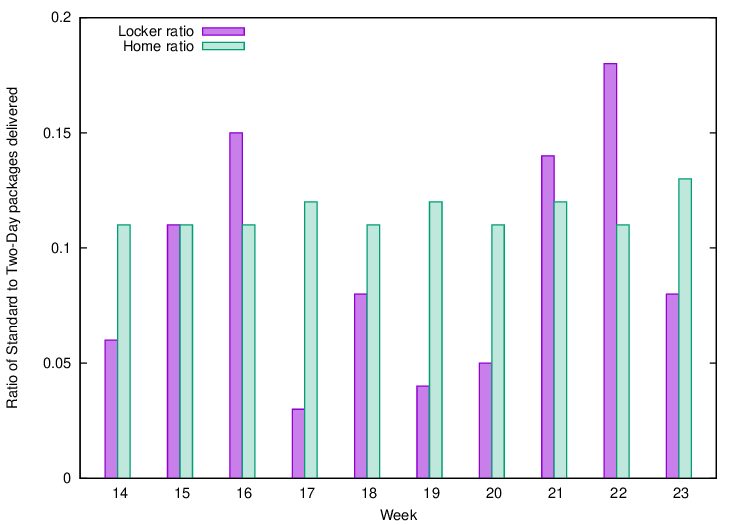}
\caption{Ratio of number of Standard to Two-Day packages delivered to 
home and locker Ruby over a ten-week period in 2018.}
\label{f_ratio_compare}
\end{figure}

This approach was an intuitive first attempt at estimating locker
demand. 
However, it had the following limitations:

\begin{itemize}
\item 
It assumes that the same ratio is maintained between the number of
packages deliveries for different ship options in homes and lockers,
which is not true in general.
For example, Figure~\ref{f_ratio_compare} shows the ratio between
number of Standard and Two-Day ship option packages delivered to homes
in zip-code 98109 and to locker Ruby over a 10-week period in 2018. 
Note that the ratio for
Locker packages has a higher variability than the ratio for
home packages. 
The two ratios can differ substantially.
\item 
Space is reserved in the ratio of the number
of packages expected to be delivered on a particular day. However,
package dwell time is not considered. The number of packages delivered
is not equal to the number of packages in a locker on a
given day. There are packages in the locker that were delivered in
the previous few days and were not picked up by the customer.
Therefore, space needs to be reserved while taking those packages
into account. Figure~\ref{f_package_dwell_time}
shows non-zero average dwell time for different ship options over
a 10-week period for lockers in zip-code 98109. 
A package that
is picked up on the same day as delivery is counted as having a
dwell time of zero.
\item 
It does not maximize locker throughput. 
It is possible
to achieve higher throughput by allocating protection limits based
on demand as well as dwell time of ship options. 
As mentioned
previously, dwell time of packages belonging to different ship
options follow a distinctive trend with higher ship option packages
having a lower dwell time (as shown in Figure~\ref{f_package_dwell_time}). 
Therefore, a model that considers demand
as well as dwell time is required to maximize throughput through
the locker.
\end{itemize}  

Amazon Locker influences customer experience directly, thus giving
rise to a high impact application with unique challenges. 
The work
presented in this paper analyzes the strengths of traditional
operations research methods and machine learning techniques, and
demonstrates best use of both worlds, thus establishing the importance to
Amazon of the rich field of yield management.

We next describe our model for reserving capacity for ship options in
Amazon Locker.

\section{Locker Capacity Management Model} 
\label{s_model}

\begin{figure}[ht]
\centering
\begin{tabular}{c}
\hline
\multicolumn{1}{|c|}{\ }\\
\multicolumn{1}{|c|}{\includegraphics[width=120 mm]{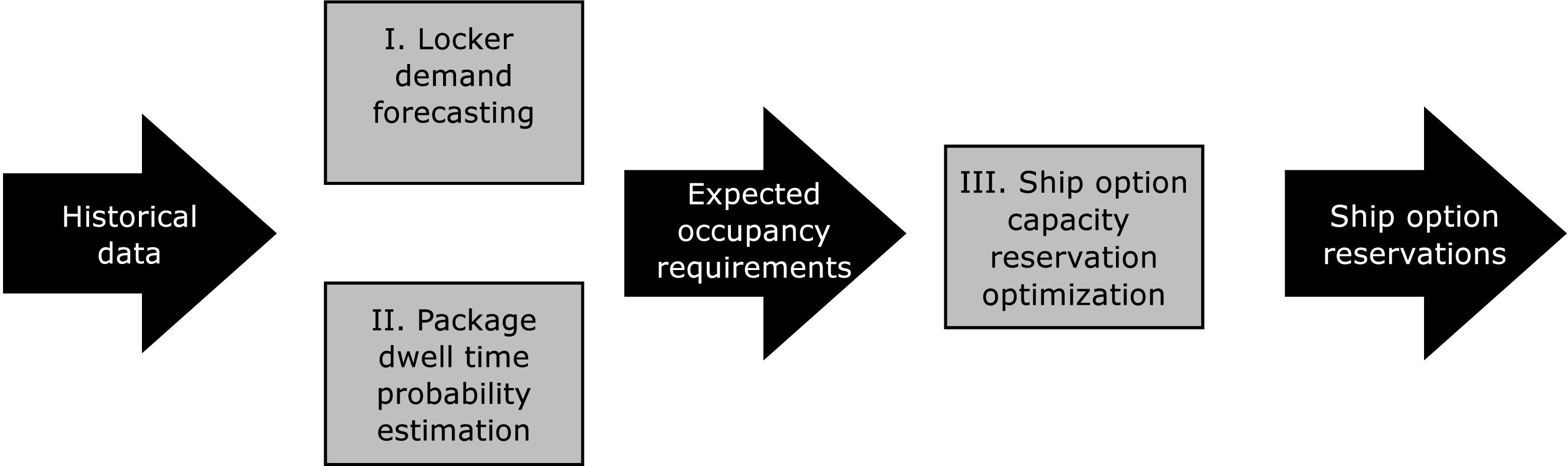}} \\
\multicolumn{1}{|c|}{\ }\\
\hline
\end{tabular}
\caption{Locker Capacity Management Model overview}
\label{f_model_overview}
\end{figure}

The locker capacity management model for assigning capacity
reservations to different ship options consists of three modules
(see Figure~\ref{f_model_overview}). 
We describe each module in detail
in Sections~\ref{ss_forecast}, \ref{ss_dwell}, and~\ref{ss_final}.

\subsection{Locker Demand Forecast} 
\label{ss_forecast}

The first module forecasts the expected number of packages delivered
in each locker, of each ship option, and on each day over the next
seven days (Module I in Figure~\ref{f_model_overview}).
\textcolor{black}{
In Section~\ref{ss_final} we refer to the outputs from this model as 
$d_{st}$, the demand forecast for ship option $s$ on day $t$.}
This problem poses some unique challenges to traditional time-series
forecasting models:
\begin{itemize}
\item 
Though data on orders that were placed and delivered are available,
data for requests that were rejected due to lack of locker capacity
are noisy and cannot be used for demand prediction. 
Therefore,
historical demand is constrained by locker capacity and previously
defined protection limits. Most traditional demand-unconstraining
methods require some portion of historical data to be unconstrained,
which is not possible in the case of lockers.
\item 
Time-series models require data for the previous few years to account
for seasonality. However, the number of lockers can increase by up
to $50\%$ year over year. Therefore, only a few months of data is
available for many lockers during peak season.
\end{itemize}

In the new method we propose, random forest regression is used to
estimate locker demand for each locker, ship option, and day, for
the next seven days. 
We train seven different random forest models, one for each day in the 
next seven days because there is more
information available for predicting demand for 
the next day than for predicting demand 
for seven days out.
The six features of the random forest regression model are:
\begin{enumerate}
\item 
Locker deliveries for that ship option on that weekday
of the previous four weeks;
\item 
Home deliveries in the Locker zip-code in the same week of the
previous year to account for seasonality in demand;
\item 
Time of first rejection for a ship option to account for unseen demand;
\item 
Delivery day of week;
\item 
Delivery day of month;
\item 
Ship option.
\end{enumerate}

The training time period is sixteen weeks of data before the model
run date, along with 
four weeks of data from peak of the previous
year.
\textcolor{black}{
This is so the model is given a broader range of possible case scenarios that
could allow it to pick up on sudden peaks that might occur within 
the planning horizon.} 
Random forest regression has the advantage of being robust
in the face of sparse data.
\textcolor{black}{
We say data is sparse since lockers are implemented in a number of
locations with different demographics (e.g., on university campuses and around
suburban neighborhoods), so the pick up and ordering characteristics
vary in volume and timing. Campuses, for example, have peak at the start
and end of the semester a time that
suburban neighborhoods may not experience peak.
} 
Furthermore, machine learning techniques
allow for additional features, such as 
\textcolor{black}{
the first observed rejection
time for unconstraining demand,
i.e., 
the first point in time when we reject a request for a delivery to the Locker 
due to lack of capacity. Since we do not have complete visibility into 
lost demand, this time directionally indicates how much demand exceeds 
capacity. An early time indicates demand is much higher than capacity while a
later time indicates that the difference between demand and 
capacity is smaller.}
Seasonality is incorporated using
home deliveries in the previous year in the same week, as one of
the features.

To evaluate model performance, prediction accuracy of the model is
compared with the values obtained with the legacy method of using
the proportion of ship option deliveries method, described in
Section~\ref{s_current_practice}.

\begin{figure}[p]
\begin{subfigure}
  \centering
  \includegraphics[width=1.0\linewidth]{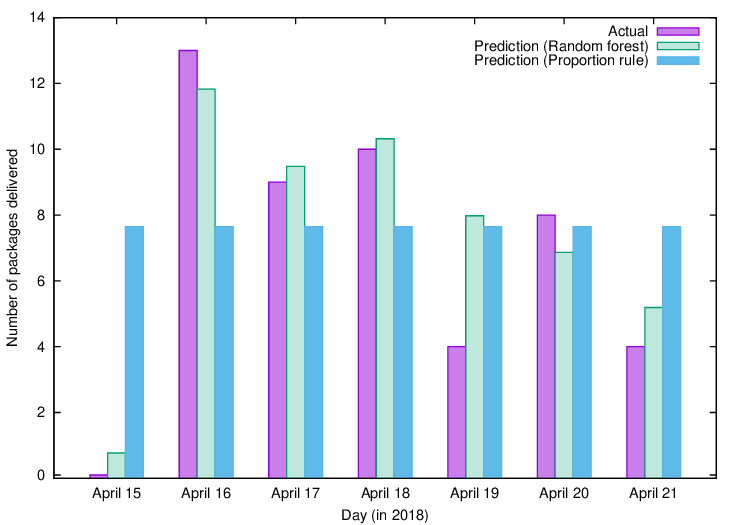}
  \label{f_ruby}
\end{subfigure}%
\begin{subfigure}
  \centering
  \includegraphics[width=1.0\linewidth]{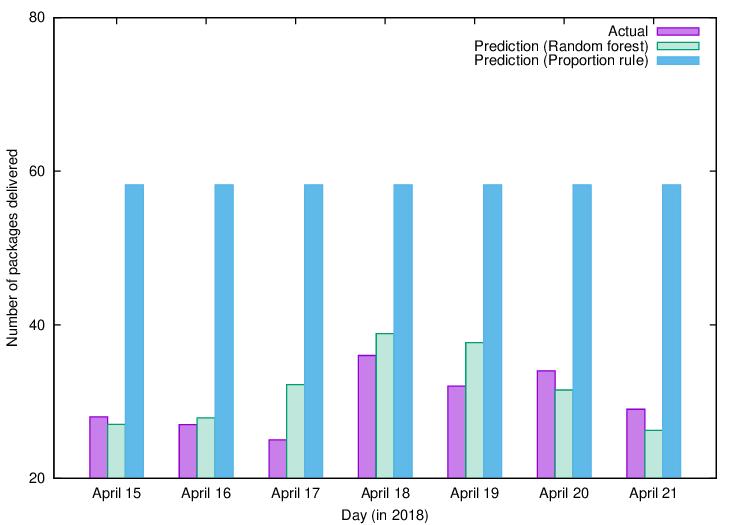}
  \label{f_lambda}
\end{subfigure}
\caption{Baseline (proportion rule) and random forest forecasts versus actual
number of packages delivered to Amazon Locker Ruby (top) and Lambda (bottom)}
\label{f_forecast}
\end{figure}

The error metric we use is not the traditional
the mean absolute percentage error (MAPE). 
This is a distinguishing feature 
of our model. 
\textcolor{black}{
The MAPE we use is defined as the average of
\begin{displaymath}
\frac{|\text{Actual Demand} - \text{Predicted Demand}|}{\text{Locker Capacity}} 
\end{displaymath}
over all lockers.
}
The reason for this is that if we predicted a demand of two packages
and there was one package 
delivered, the regular MAPE would say that there was a 100\% error. 
However, while this error is significant if locker capacity is three slots, 
it is less significant if Locker capacity is 100 slots. 
The new error metric \citep{Bruce2019} uses
the same units for model performance
as the business performance metric,
locker throughput.
For a two-week test period (April 15 to 28, 2018), 
the MAPE 
decreased from $12.4\%$ for the
\textcolor{black}{legacy}
method to $2.9\%$ for the random forest regression, for a test set
of $30$ lockers in the US. MAPE values for individual lockers, and
two ship options are presented in Figure~\ref{f_forecast}, where
plots of prediction versus actual for both methods, for two
sample lockers, are shown. 
The disruptive improvement obtained using random
forest regression is clear.

\subsection{Package Dwell Time Probability Estimation} \label{ss_dwell}

Module~II of the Locker Capacity Management Model
(see Figure~\ref{f_model_overview}) is
package dwell time probability estimation.
In this module, we estimate the probability that a package stays
in the locker for $0,1,2, \ldots, 6$ days for each shipping speed
and delivery day in the time horizon of the next seven days. 
\textcolor{black}{
In Section~\ref{ss_final},
we refer to the output dwell time probabilities from this model as $p_{svt}$, 
the probability that a package belonging to ship option $s$ and delivered 
on day $v$ will be present on day $t$, where $t > v$.}

Customers are advised to pick up their packages within three days
of delivery. 
If the package is not picked up in three days, a call
tag is generated and a carrier picks up the package for return to
Amazon within the next two to three days. 
Therefore, the possible values for
dwell time are $0,1,\ldots,6$. 
Dwell time probability is currently
calculated in the locker capacity management system for a separate
function unrelated to the capacity reservations for ship options.
This probability is calculated as the proportion of packages that
stayed in the locker for $0,1,2, \ldots, 6$ days in the previous
four weeks. 
However, this method often leads to over-fitting due
to sparse data. 
For example, if in the previous four weeks, the
number of Standard deliveries is two, and one package was picked
up on the same day as delivery, and the other stayed in the locker
for six days, with the prior method, the probability of the package
being picked up on day $0$ will be $0.5$, and the probability the
package stays in the locker for six days will be $0.5$, while the
probability that the package stays for $1, 2,3,4$ and $5$ days is
$0$, which is not descriptive of reality. 
In reality, the probability
that a package stays in the locker for one day is not zero, and is
higher than the probability the package stays in the locker for
six days.

In this module, random forest classification is used to 
estimate dwell time probabilities of a package with the following four 
features:
\begin{enumerate}
\item Average, minimum, and maximum dwell time of packages for 
that ship option, delivered on that day of week in the previous four weeks;
\item Ship option;
\item Delivery day of week;
\item Delivery day of month.
\end{enumerate}

Random forest regression
allows utilization of data for multiple similar lockers 
to estimate dwell time probability of a package, 
thus reducing effects of sparse data.
To improve the quality of dwell time probabilities, we utilize
probability calibration using isotonic regression implemented in
the scikit-learn Python library with the technique proposed by
\citet{zadrozny2002transforming}.

\textcolor{black}{
Isotonic regression is a type of regression analysis used for fitting
a non-decreasing function to data. In the context of machine learning,
isotonic regression is often used to calibrate the probabilities
output by a classification model, making them more interpretable
and reliable.
}

\textcolor{black}{
The scikit-learn Python library provides an implementation of
isotonic regression that can be used for various applications. The
technique proposed by Zadrozny and Elkan in 2002 is specifically
geared towards probability calibration. In their method, they aim
to adjust the predicted probabilities such that they better represent
the true probabilities. The idea is to fit a non-decreasing function
to the output of a model, in a way that minimizes some measure of
error with respect to the true, observed outcomes. The implementation
in scikit-learn often follows a ``pool adjacent violators" algorithm
for this fitting process.
}

\textcolor{black}{
Using isotonic regression for calibration makes sense when you want
the output probabilities of your model to be well-calibrated, meaning
that a predicted probability of $x$\% should correspond to the event
happening $x$\% of the time. This is especially useful in risk-sensitive
applications where the predicted probabilities are used to make
decisions.
}

\textcolor{black}{
To use isotonic regression in scikit-learn, you typically first
train your classification model and obtain predicted probabilities.
Then, you fit an isotonic regression model to these predicted
probabilities and the true labels, thereby obtaining a calibrated
set of probabilities.
}


Similar to the error metric in Section~\ref{ss_forecast},
the error metric used
to measure the quality of probabilities also uses
the same units
in its model performance metric as used in the business performance metric,
i.e. locker throughput. 
The performance metric is 
\textcolor{black}{
the average of
\begin{displaymath}
\frac{|\text{Actual \# of Packages Picked Up} - \text{Expected \# of Packages Picked Up}|}{\text{Locker Capacity}}, 
\end{displaymath}
over all lockers and over the days in the planning time horizon of
where the expected number of packages picked up is calculated using
the probabilities computed in this module and the actual number
of packages delivered to the locker. 
}
When compared with a traditional
metric such as log loss, this performance metric has the advantage
of translating improvements to the model performance directly to
improvements in locker utilization predictions and makes it easier
to explain the model performance to business teams. 
Compared
with the methodology of assuming that all packages would be picked
up on the same day as delivery, this method resulted in an eight
percent overall improvement in the error metric for lockers in the US in a
two-week test period (April 15 to 28, 2018).

\subsection{Ship Option Capacity Reservation Optimization} 
\label{ss_final}

In the ship option capacity reservation module (Module~III 
of Figure~\ref{f_model_overview}), 
a linear programming
formulation is proposed to find optimal capacity reservations for
ship options in a locker with the goal of maximizing throughput. 
This linear program
is run at the end of day $0$, and the time horizon is assumed to
be the next seven days, with the count starting with day $1$. 
Based on
the customer or carrier pickup time windows described
in 
Section~\ref{ss_dwell}, it 
is assumed that packages delivered
up to six days before the first day of the time horizon might still
be in the locker at the end of day $0$. 
The actual number of packages
in the locker is known in real time and is a parameter. 
However,
the date on which the package will be picked up is not known.

Let $T$ be the number of days in the time horizon and $S$ be the number of 
ship options.
Let $C$ be the locker capacity and $p_{svt}$ be the
\textcolor{black}{dwell time} probability of a package delivered for 
ship option $s \in \{1, 2, \ldots, S \}$
on day $v \in \{-6,-5, \ldots, t\}$ being present in the locker on day 
$t \in \{1, 2, \ldots, T\}$. 
This probability is computed in Module~II of
Figure~\ref{f_model_overview}. 
Furthermore, let
$e_{sv}$ \textcolor{black}{be the number of packages} of ship option $s \in \{1,2, \ldots, S\}$ 
delivered to the locker on day $v \in \{-6,-5,\ldots,0\}$ and present in the 
locker at
the end of day 0 and
$d_{st}$ is the \textcolor{black}{package demand forcast} for package deliveries of ship option $s$ on day 
$t \in \{1, \ldots, T\}$.
This demand is computed in Module~I of 
Figure~\ref{f_model_overview}, 

There are two classes of decision variables.
Decision variable $x_{st}$ is the  number of slots to reserve for ship option 
$s \in \{1,2,\ldots,S\}$ on day $t \in \{1, 2, \ldots, T\}$.
Variable $y_{st}$ is the number of packages accepted for arrival for 
ship option $s \in \{1, 2, \ldots, S\}$  on day 
$t \in \{1, 2, \ldots, T\}$.

A linear program whose objective is to maximize throughput,
i.e. the total number of packages delivered to the locker in the
time horizon, is given in Equations~(\ref{eq1a})--(\ref{eq1e})
in the Appendix.

In the linear program,
decision variables $x_{st}$ are used as capacity reservations for 
ship option $s \in \{1, 2, \ldots, S\}$ on day 
$t \in \{1, \ldots, T\}$ of the
of planning horizon. 

In this formulation, faster ship option packages are not
explicitly prioritized
using specific weight parameters. 
However, since the average dwell
times of faster ship option packages is lower than that of slower
ship option packages (see Figure~\ref{f_package_dwell_time}),
the objective of throughput maximization automatically results in
prioritization of faster ship option packages.

\begin{figure}[ht]
\centering
\includegraphics[width=120 mm]{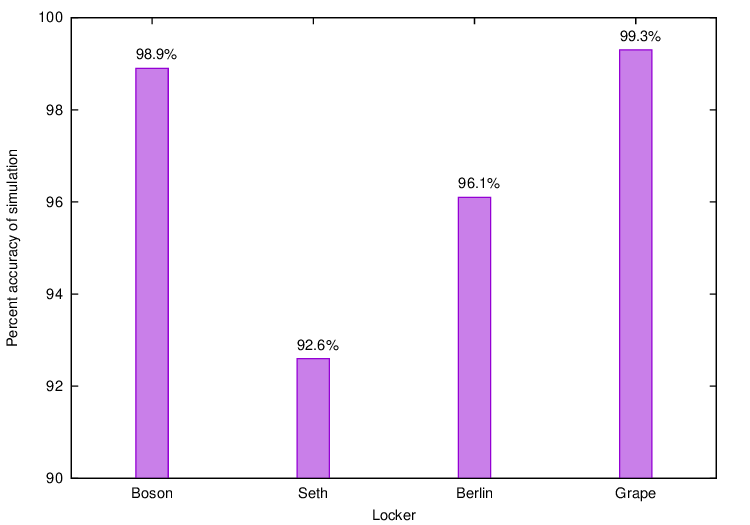}
\caption{Simulation system accuracy metrics}
\label{f_simulation}
\end{figure}

\begin{figure}[ht]
\centering
\includegraphics[width=100mm]{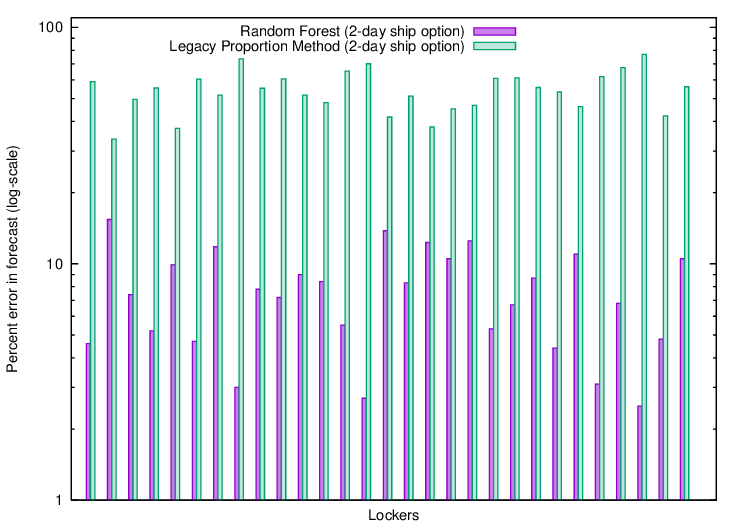}
\caption{Percent error in package demand forecast for random forest method
and legacy proportion 
method on two-day ship option.}
\label{f_accuracy-30-lockers-a}
\end{figure}

\begin{figure}[ht]
\centering
\includegraphics[width=100mm]{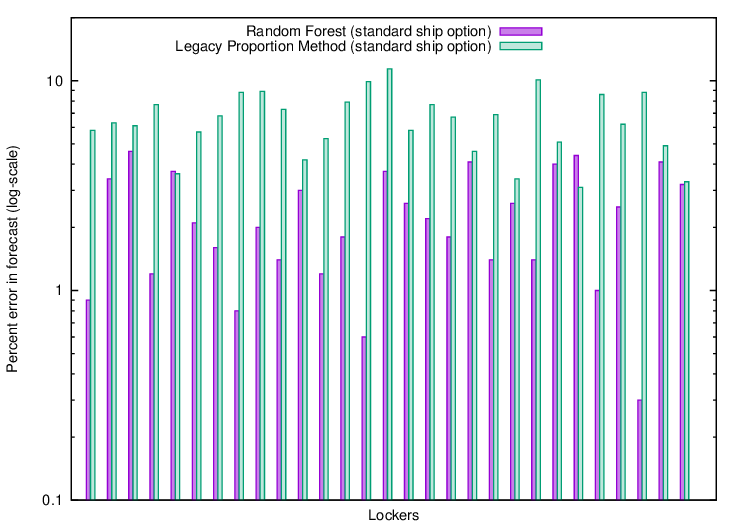}
\caption{Percent error in package demand forecast for random forest method
and legacy proportion 
method on standard ship option.}
\label{f_accuracy-30-lockers-b}
\end{figure}

\begin{figure}[ht]
\centering
\includegraphics[width=120 mm]{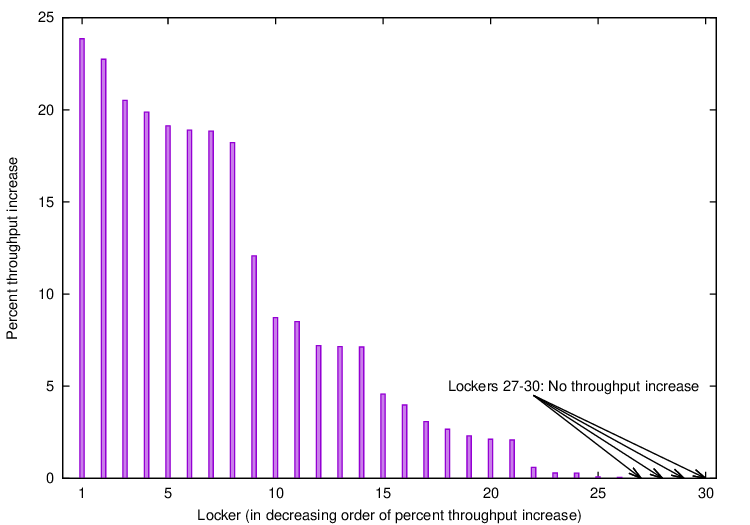}
\caption{Percent throughput increase with new Locker Capacity Management
algorithm compared to legacy rule.
The 30 lockers are sorted in decreasing order of throughput increase.}
\label{f_throughput-increase}
\end{figure}

\section{Model Evaluation by Simulation System}
\label{s_simulation}

To best serve Amazon customers,
it is essential to quantify the
impact of any improvement to the capacity management system. 
However,
locker capacity management is a complex system, and it is difficult
to predict the effect of any change to inputs on locker throughput.
To enable this task,
a simulation system was deployed allowing the
determination of the
impact
on the locker
capacity management system
of any changes made to inputs, 
including changes in ship option capacity reservations.

This simulation system was built to mirror what happens in production,
based on different inputs. 
Millions of package events (e.g. order requests,
arrivals, order pick ups, etc.) that happened  in a locker
over a 15-day period can be replayed
in this simulation system in less than six minutes. 
This tool has enabled us to
perform many types of \textit{what-if} analyses.
Prior
to building the simulation system, the only way to do these analyses
was to apply changes directly to production and measure the result in
real-life. 
This is now
done through the simulation system, thus
reducing the time needed to arrive at data-driven decisions, and
eliminating the risk of negatively impacting customer experience.
Accuracy of the simulation system is measured as the percentage of
acceptance/rejection decisions of locker requests that were the
same in production and in the simulation system. 
Figure~\ref{f_simulation}
shows the high level of simulation system accuracy for 
four sample lockers over a
$14$-day time period.

\section{Experimental Results} 
\label{s_results}

In this section we report
experimental results from the implementation of the 
Locker Capacity Management Model
both in the simulation system and in \textit{real-world} production.
Results are presented in terms of change in
locker throughput,
the adopted
business metric for measuring the success of the new model.

\subsection{Simulations} \label{ss_sim_results}

The ship option capacity reservations from the model proposed in
Section~\ref{s_model} were first tested on the simulation system
for the time period of the two weeks from April 15 to
28, 2018. 
A set of 30 lockers were used in the experiment, 
chosen randomly as a mix of high, medium, and
low throughput lockers in the US.
\textcolor{black}{
Note that we used actual data about the lockers as well as historical data. 
}
Two sets of simulations were run for each
locker. 
The first set of simulations used the ship option capacity reservations 
previously deployed
in production,
\textcolor{black}{
i.e. the legacy system}.
The second set of simulations used the ship option capacity reservations
obtained through the algorithm described in this paper. 
Historical
locker requests for the 30 lockers were replayed, and the number of
accepted requests in the two sets of simulations were compared
against each other. 
Simulation results are shown in 
Figures~\ref{f_accuracy-30-lockers-a} to~\ref{f_throughput-increase}.
Figure~\ref{f_accuracy-30-lockers-a} and~\ref{f_accuracy-30-lockers-b} 
compare, respectively,
the percent error (in log scale) in the forecasts made with random 
forest regression and with the legacy proportion method on the
lockers used in the study for two-day ship option 
and standard ship option.
The percent error are computed as the absolute difference between the 
forecast and what actually occurred divided by the locker capacity.
For example, if the forecast was 1 and the actual was 2 while the locker
capacity was 100, then the percent error would be $|2-1|/100 = 1$ percent.
However, if the locker capacity were 3, then the percent error would be 
$|2-1|/3 =33.33$
percent.  

Note that while random forest regression forecasts are better than those
made with the legacy approach, there is a more noticeable difference for 
the two-day ship option.
On the standard ship option both methods do well with random forests still
outperforming the legacy method.

Figure~\ref{f_throughput-increase} shows percent throughput increase with the
new capacity management algorithm as compared to the legacy algorithm.
There is an average $6\%$ increase in
throughput for the 30 lockers with a maximum increase of $23\%$.
In four of the 30 lockers there was no increase or decrease in locker 
throughput.
This is because those lockers had low demand when compared to
locker capacity, so any model, including first-come-first-served
would result in the same throughput.
These simulations provided the needed confidence in the model
to launch on actual lockers in production.

\subsection{Lockers in Production} \label{ss_results}

\begin{figure}[ht]
\centering
\includegraphics[scale=1]{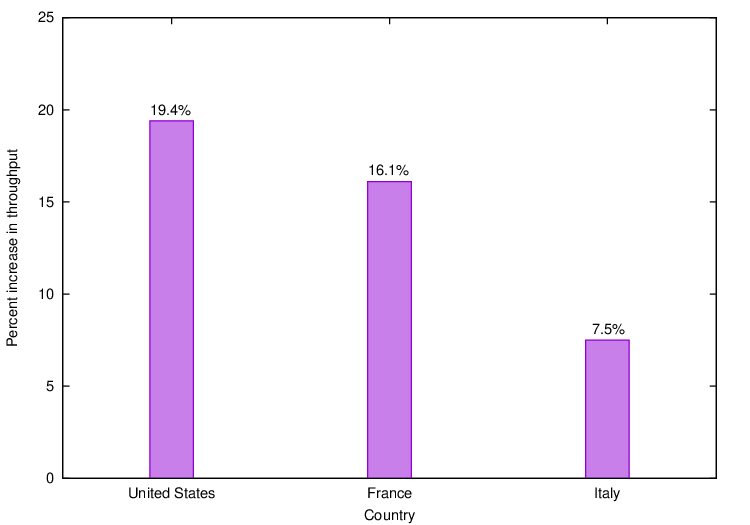}
\caption{Model implementation results show for the US, FR and IT the percent 
increase in locker throughput achieved by deploying the new Locker Capacity
Management algorithm compared to the legacy rule.}
\label{f_model_implementation}
\end{figure}

The model was put into production on live lockers in a phased deployment,
with 100\% of Amazon Lockers worldwide
utilizing model outputs by November 2018. 
This resulted in an increase in throughput in
all countries. 
Note, however, that if a locker has very low demand, 
then our model is not expected to make much of a difference, since 
the first-come-first-served policy performs as well as any model. 
Therefore the
results were evaluated on Amazon Locker locations having 
\textcolor{black}{demand to
capacity ratio} above the 
60 percentile
mark.
Here, these percentile cutoffs are calculated based on all Amazon 
Lockers in the country for the
five weeks during peak.
Figure~\ref{f_model_implementation} shows percent improvement results in the
US, France, and Italy.
This allowed Amazon to say ``yes'' to
millions of requests for which it would otherwise reject.

\section{Conclusion} 
\label{s_conclusion}

Amazon Locker influences customer experience directly, thus giving
rise to a high impact yield management application with unique
challenges due to data sparsity and unknown package dwell times.
The major pain point for Amazon Locker customers is rejected requests due
to capacity algorithm estimations. The biggest chunk ($40\%$) of
unjustified rejections of requests for deliveries to lockers in
2017 were due to poor ship option capacity reservations. Amazon
built and implemented an algorithm that computes optimal ship option
capacity reservations to maximize throughput (i.e. minimize rejections)
and therefore, enhance customer experience at Amazon Lockers. 
This allowed
Amazon to accept and say ``yes'' to millions of additional Locker
requests.

\textcolor{black}{
The methodology introduced in this paper and utilized by us for managing 
parcel locker capacity has
versatile applications across various business scenarios, particularly
in optimizing space and resource allocation. 
For instance, in a
grocery store setting, this approach could enhance
the management of shelf space, especially for perishable goods where
shelf life and demand prediction are crucial. The model can dynamically
adjust inventory levels and shelf allocation based on real-time
demand forecasts, reducing waste and maximizing efficiency. Extending
this idea further, such innovative methods can revolutionize the
management of computational resources in cloud services like AWS.
Here, the approach could be tailored to allocate CPU services more
effectively. Customers often request compute time without precise
knowledge of the duration needed for their computations. By employing
a system that intelligently bundles and allocates CPU resources
across the cloud based on demand and cost considerations, efficiency
and cost-effectiveness can be significantly improved. 
}

\textcolor{black}{
One significant limitation inherent in these systems is their
dependency on accurate forecasting and an understanding of customer
behavior. This reliance becomes particularly challenging due to the
dynamic and often unpredictable nature of consumer preferences and
market trends. If given more time and resources, a
more nuanced model could be developed. This model would not only
factor in general buying trends but also delve into customer behavior
specific to each product (identified by its ASIN or SKU). Such a granular
approach could potentially yield more precise predictions and
efficient allocation of resources.
}

\textcolor{black}{
Furthermore, an intriguing avenue for exploration could be the
strategic overbooking of space, akin to practices in the airline
and hotel industries. This approach could potentially maximize
utilization and revenue, especially in scenarios where demand
outstrips supply. However, it is crucial to balance this strategy
with the risk of overcommitment, which could lead to customer
dissatisfaction and operational challenges. Assessing the threshold
for overbooking and developing mechanisms to manage potential
overcommitment situations are key areas for further research
and development.
}

\textcolor{black}{
Considering future market trends, we have identified 
the interplay between locker deliveries and last-mile
delivery services as
a
promising area for exploration that has not yet been thoroughly
investigated. 
Understanding customer preferences and behavior
in this domain could yield valuable insights. Specifically, it would
be intriguing to conduct a study where customers are given a choice
between immediate locker delivery and next-day home delivery. This
investigation could reveal customer preferences for convenience versus
immediacy, and how these choices align with their expectations and
satisfaction levels.
}

\textcolor{black}{
Moreover, the financial aspect of this comparison is equally
important. We hypothesize that home deliveries, while offering
direct-to-door convenience, likely incur higher costs due to factors
such as fuel, labor, and vehicle maintenance. In contrast, locker
deliveries, by centralizing drop-offs, could significantly reduce
these expenses. Assessing the cost implications of each delivery
method would not only inform our operational strategies but also
help in tailoring our services to be both cost-effective and
customer-centric.
}

\textcolor{black}{
The potential environmental impact
of each delivery method
could also be explored. 
For instance, locker deliveries might
reduce carbon emissions by reducing the distance covered by
delivery vehicles. This aspect aligns with growing consumer awareness
and concern about environmental sustainability.
}

\section*{Acknowledgements}

The work of Sethuraman, Mardan, and Resende was mainly done when they were
at Amazon Transportation Services in Seattle and Bellevue, WA.
The work of Bansal was mainly done when he was at Access Point Technology, 
Amazon Delivery Technologies, in Delhi, India.

On December 13, 2022 a United States patent was issued for the
capacity management system described in this paper \citep{patent2022}.

\bibliographystyle{plainnat}
\bibliography{references}

\begin{thebibliography}{29}
\providecommand{\natexlab}[1]{#1}
\providecommand{\url}[1]{\texttt{#1}}
\expandafter\ifx\csname urlstyle\endcsname\relax
  \providecommand{\doi}[1]{doi: #1}\else
  \providecommand{\doi}{doi: \begingroup \urlstyle{rm}\Url}\fi

\bibitem[Beck et~al.(2021)Beck, Harvey, Kaylen, Sala, Urban, Vermeulen, Wilken,
  Xie, Iliescu, and Mital]{BecEtAlii21a}
J.~Beck, J.~Harvey, K.~Kaylen, C.~Sala, M.~Urban, P.~Vermeulen, N.~Wilken,
  W.~Xie, D.~Iliescu, and P.~Mital.
\newblock Carnival optimizes revenue and inventory across heterogenous cruise
  line brands.
\newblock \emph{INFORMS Journal on Applied Analytics}, 51\penalty0
  (1):\penalty0 26--41, 2021.

\bibitem[Bitran and Caldentey(2003)]{BitCal03a}
G.~Bitran and R.~Caldentey.
\newblock An overview of pricing models for revenue management.
\newblock \emph{Manufacturing \& Service Operations Management}, 5:\penalty0
  203--229, 2003.

\bibitem[Bruce(2019)]{Bruce2019}
A.~Bruce.
\newblock Personal communication, 2019.

\bibitem[Clough et~al.(2015)Clough, Jacobs, and Gel]{CloJacGel15a}
{M.C.} Clough, {T.L.} Jacobs, and E.~Gel.
\newblock New formulations for price and ticket availability decisions in
  choice-based network revenue management.
\newblock In \emph{{AGIFORS 55th Annual Symposium: Analytics for Efficiency and
  Customer Centric Optimization}}, Washington, D.C., August 2015. Airline Group
  of the International Federation of Operational Research Societies (AGIFORS).

\bibitem[Feldman and Topaloglu(2017)]{FelTop17a}
J.B. Feldman and H.~Topaloglu.
\newblock Revenue management under the {Markov Chain} choice model.
\newblock \emph{Operations Research}, 65:\penalty0 1322--1342, 2017.

\bibitem[Ferreira et~al.(2018)Ferreira, {Simchi-Levi}, and Wang]{FerSimWan18a}
K.J. Ferreira, D.~{Simchi-Levi}, and He~Wang.
\newblock Online network revenue management using thompson sampling.
\newblock \emph{Operations Research}, 66:\penalty0 1586--1602, 2018.

\bibitem[Gopalakrishnan and Rangaraj(2010)]{GopRan10a}
R.~Gopalakrishnan and N.~Rangaraj.
\newblock Capacity management on long-distance passenger trains of {Indian
  Railways}.
\newblock \emph{Interfaces}, 40\penalty0 (4):\penalty0 291--302, 2010.

\bibitem[Gorman(2010)]{Gor10a}
M.F. Gorman.
\newblock {Hub Group implements a suite of OR tools to improve its operations}.
\newblock \emph{Interfaces}, 40\penalty0 (5):\penalty0 368--384, 2010.

\bibitem[Guillen et~al.(2019)Guillen, Ruiz, Dellepiane, Maccarrone, Maccioni,
  Pinzuti, and Procacci]{GuiRuiDelMacMacPinPro19a}
J.~Guillen, P.~Ruiz, U.~Dellepiane, L.~Maccarrone, R.~Maccioni, A.~Pinzuti, and
  E.~Procacci.
\newblock Europcar integrates forecasting, simulation, and optimization
  techniques in a capacity and revenue management system.
\newblock \emph{INFORMS Journal on Applied Analytics}, 49\penalty0
  (1):\penalty0 40--51, 2019.

\bibitem[Iwan et~al.(2016)Iwan, Kijewska, and Lemke]{IwaKijLem16a}
S.~Iwan, K.~Kijewska, and J.~Lemke.
\newblock {Analysis of Parcel Lockers' Efficiency as the Last Mile Delivery
  Solution -- The Results of the Research in Poland}.
\newblock In \emph{Tenth International Conference on City Logistics 17-19 June
  2015}, volume~12 of \emph{Transportation Research Procedia}, pages 644--655,
  2016.

\bibitem[Jacobs et~al.(2008)Jacobs, Smith, and
  Johnson]{jacobs2008incorporating}
T.L. Jacobs, B.C. Smith, and E.L. Johnson.
\newblock Incorporating network flow effects into the airline fleet assignment
  process.
\newblock \emph{Transportation Science}, 42\penalty0 (4):\penalty0 514--529,
  2008.

\bibitem[Kunnumkal and Talluri(2019)]{KunTal19a}
S.~Kunnumkal and K.~Talluri.
\newblock A strong {Lagrangian} relaxation for general discrete-choice network
  revenue management.
\newblock \emph{Computational Optimization and Applications}, 73:\penalty0
  275--310, 2019.

\bibitem[Littlewood(1972)]{Lit72a}
K.~Littlewood.
\newblock Forecasting and control of passenger bookings.
\newblock In \emph{Proceedings of 12th AGIFORS Symposium}, pages 95--117,
  Nathanya, Israel, October 1972.

\bibitem[Littlewood(2005)]{littlewood}
K.~Littlewood.
\newblock Forecasting and control of passenger bookings.
\newblock \emph{Journal of Revenue and Pricing Management}, 4\penalty0
  (2):\penalty0 111--123, 2005.

\bibitem[McGill(1989)]{McG89a}
J.I. McGill.
\newblock \emph{Optimization and estimation problems in airline yield
  management}.
\newblock PhD thesis, University of British Columbia, December 1989.

\bibitem[McGill and {van Ryzin}(1999)]{McGRyz99a}
J.I. McGill and G.J. {van Ryzin}.
\newblock Revenue management: {Research} overview and prospects.
\newblock \emph{Transportation Science}, 33:\penalty0 233--256, 1999.

\bibitem[Orenstein et~al.(2019)Orenstein, Raviv, and Sadan]{OreRavSad19a}
I.~Orenstein, T.~Raviv, and E.~Sadan.
\newblock {Flexible parcel delivery to automated parcel lockers: Models,
  solution methods and analysis}.
\newblock \emph{EURO Journal on Transportation and Logistics}, 8:\penalty0
  683--711, 2019.

\bibitem[Pandey et~al.(2017)Pandey, Dutta, and Joshi]{PanDutJos17a}
S.~Pandey, G.~Dutta, and H.~Joshi.
\newblock Survey on revenue management in media and broadcasting.
\newblock \emph{Interfaces}, 47\penalty0 (3):\penalty0 195--213, 2017.

\bibitem[Phillips(2005)]{phillips2005pricing}
R.L. Phillips.
\newblock \emph{Pricing and revenue optimization}.
\newblock Stanford University Press, 2005.

\bibitem[Rohmer and Gendron(2020)]{RohGen20a}
S.~Rohmer and B.~Gendron.
\newblock {A guide to parcel lockers in last mile distribution -- Highlighting
  challenges and opportunities from an OR perspective}.
\newblock Technical Report CIRRELT-2020-11, CIRRELT, Universit\'e de
  Montr\'eal, 2020.

\bibitem[Sethuraman et~al.(U.S. Patent 11 526 838, Dec. 2022)Sethuraman,
  Mardan, Jacobs, Jain, and Bansal]{patent2022}
S.~Sethuraman, S.~Mardan, T.~L. Jacobs, M.~Jain, and A.~Bansal.
\newblock Capacity management system for delivery lockers, U.S. Patent 11 526
  838, Dec. 2022.

\bibitem[Shihab et~al.(2019)Shihab, Logemann, Thomas, and Wei]{ShiLogThoWei19a}
S.A.M. Shihab, C.~Logemann, D.-G. Thomas, and Peng Wei.
\newblock Autonomous airline revenue management: {A} deep reinforcement
  learning approach to seat inventory control and overbooking.
\newblock In \emph{Reinforcement Learning for Real Life ({RL4RealLife})
  Workshop in the 36th {International Conference on Machine Learning}}, Long
  Beach, California, 2019.

\bibitem[Smith et~al.(1992)Smith, Leimkuhler, and Darrow]{smith1992yield}
B.C. Smith, J.F. Leimkuhler, and R.M. Darrow.
\newblock Yield management at {American Airlines}.
\newblock \emph{Interfaces}, 22\penalty0 (1):\penalty0 8--31, 1992.

\bibitem[Talluri et~al.(2009)Talluri, {van Ryzin}, Karaesmen, and
  Vulcano]{TalRuzKarVul09a}
K.T. Talluri, G.J. {van Ryzin}, I.Z. Karaesmen, and G.J. Vulcano.
\newblock Revenue management: {Models} and methods.
\newblock In \emph{Proceedings of the 2009 Winter Simulation Conference}, pages
  148--161, Austin, Texas, 2009.

\bibitem[{van Duin} et~al.(2020){van Duin}, Wiegmans, {van Arem}, and {van
  Amstel}]{DuiWieAreAms20a}
J.H.R. {van Duin}, B.W. Wiegmans, B.~{van Arem}, and Y.~{van Amstel}.
\newblock From home delivery to parcel lockers: {A} case study in amsterdam.
\newblock \emph{Transportation Research Procedia}, 46:\penalty0 37--44, 2020.

\bibitem[Wang(1983)]{wang1983optimum}
K.W. Wang.
\newblock Optimum seat allocation for multi-leg flights with multiple fare
  types.
\newblock In \emph{Proceedings of the Twenty-third annual symposium - AGIFORS,
  Airline Group of the International Federation of Operational Research
  Societies}, Olive Branch, Mississippi, 1983.

\bibitem[Wang et~al.(2018)Wang, Liu, Chen, Ren, Meng, and
  Hu]{WanLiuCheRenMenHu18a}
Lulu Wang, Huahui Liu, Guanhao Chen, Shaola Ren, Xiaonan Meng, and Yi~Hu.
\newblock Learning theory and algorithms for revenue management in sponsored
  search.
\newblock Technical Report 1807.01827, arXiv, 2018.

\bibitem[Weatherford and Kimes(2003)]{WeaKim03a}
L.R. Weatherford and S.E. Kimes.
\newblock A comparison of forecasting methods for hotel revenue management.
\newblock Technical report, Cornell University, School of Hotel Administration
  site, 2003.
\newblock URL \url{http://scholarship.sha.cornell.edu/articles/753}.
\newblock Retrieved from institution site on 2020-04-24.

\bibitem[Zadrozny and Elkan(2002)]{zadrozny2002transforming}
B.~Zadrozny and C.~Elkan.
\newblock Transforming classifier scores into accurate multiclass probability
  estimates.
\newblock In \emph{Proceedings of the Eighth ACM SIGKDD International
  Conference on Knowledge Discovery and Data Mining}, pages 694--699. ACM,
  2002.

\end{thebibliography}

\pagebreak
\appendix
\section*{\textcolor{black}{Appendix: LP Model for Throughput Maximimzation}}

\textcolor{black}{
This appendix describes a linear program whose objective is 
to maximize throughput,
i.e. the total number of packages delivered to the locker in the
time horizon.
Let $S$ be the number of ship options and $T$ the number of days in the
planning horizon.
}

Let $C$ be the locker capacity and $p_{svt}$ be the
\textcolor{black}{dwell time} probability of a package delivered for 
ship option $s \in \{1, 2, \ldots, S \}$
on day $v \in \{-6,-5, \ldots, t\}$ being present in the locker on day 
$t \in \{1, 2, \ldots, T\}$. 
Furthermore, let
$e_{sv}$ \textcolor{black}{be the number of packages} of ship option $s \in \{1,2, \ldots, S\}$ 
delivered to the locker on day $v \in \{-6,-5,\ldots,0\}$ and present in the 
locker at
the end of day 0 and
$d_{st}$ is the \textcolor{black}{package demand forcast} for package deliveries of ship option $s$ on day 
$t \in \{1, \ldots, T\}$.

\textcolor{black}{
In this model, there are two classes of decision variables.
Decision variable $x_{st}$ is the  number of slots to reserve for ship option 
$s \in \{1,2,\ldots,S\}$ on day $t \in \{1, 2, \ldots, T\}$.
Variable $y_{st}$ is the number of packages accepted for arrival for 
ship option $s \in \{1, 2, \ldots, S\}$  on day 
$t \in \{1, 2, \ldots, T\}$.
}

\textcolor{black}{
\begin{align}
\text{maximize} &\qquad \sum_{t=1}^T \sum_{s=1}^S y_{st}\label{eq1a}\\
\text{subject to}&\qquad \sum_{v=1}^t \sum_{s=1}^S p_{svt}y_{sv} + 
\sum_{v=-6}^0 \sum_{s=1}^S p_{svt}e_{sv} \leq C, \; \forall t \in \{1, 2, \ldots, T\} \label{eq1b}   \\
&\qquad y_{st} \leq d_{st}, \; \forall s \in \{1, 2, \ldots, S\}, \forall t \in \{1, 2, \ldots, T\} \label{eq1c}   \\
&\qquad \sum_{v=1}^t p_{svt}y_{sv} + \sum_{v=-6}^0 p_{svt}e_{sv} = x_{st}, \; \forall s \in \{1, 2, \ldots, S\}, \forall t \in \{1, 2, \ldots, T\} \label{eq1d} \\
&\qquad x_{st} \geq 0; \;  y_{st} \geq 0, \; \forall s \in \{1, 2, \ldots, S\}, \forall t \in \{1, 2, \ldots, T\}. \label{eq1e}
\end{align}
}

\textcolor{black}{
Constraints~(\ref{eq1b}) guarantee that locker capacity is not violated in
any period of the planning horizon.
Constraints~(\ref{eq1c}) ensure, for any period of the planning horizon,
that the number of packages of any ship option
accepted for delivery is less than or equal to demand for that ship option. 
Constraints~(\ref{eq1d}) define the relationship between the number
of packages accepted for delivery and the number of slots reserved
for the ship option.
Finally, the decision variables are constrained to be nonnegative by
Constaints~(\ref{eq1e}).
}

\end{document}